\documentclass[12pt]{amsart}

\usepackage{amsmath}
\usepackage{latexsym}
\usepackage{amssymb}
\usepackage{amsfonts}
\usepackage{cite}

\usepackage{graphics}
\include{PDF}

\renewcommand{\proof}{\par\noindent{\it Proof.\ \ }}
\def\qed{\ifmmode\square\else\nolinebreak\hfill
$\Box$\fi\par\vskip12pt}

\def\l{\langle} \def\r{\rangle} 
\def\div{\,\big|\,} 

    \def\ZZ{\mathbb Z}
 \def\CC{{\mathcal C}}

\def\mod{{\sf mod~}} \def\val{{\sf val}}

\def\soc{{\sf soc}}

\def\C{{\bf C}} \def\O{{\bf O}}

\def\Ga{{\it \Gamma}}

\def\Cay{{\sf Cay}}  
\def\Aut{{\sf Aut}}  \def\Out{{\sf Out}}

\def\D{{\sf D}} 
\def\S{{\sf S}} 
\def\J{{\sf J}} \def\M{{\sf M}}

\def\Row{{\sf Row}}
\def\Transitivity{{\sf Transitivity}}\def\Yes{{\sf Yes}}\def\No{{\sf No}}
\def\transitive{{\sf transitive}}

\def\Bipartite{{\sf Bipartite}}

\def\a{\alpha} \def\b{\beta} \def\g{\gamma}

\def\PSp{{\sf PSp}}
\def\GammaL{{\sf \Gamma L}}
\def\AGammaL{{\sf A\Gamma L}}
 
\def\A{{\sf A}}

\def\PSL{{\sf PSL}}\def\PGL{{\sf PGL}}
 \def\SL{{\sf SL}}

\def\ASL{{\sf ASL}}
\def\AGL{{\sf AGL}}
\def\ASigmaL{{\sf A\Sigma L}}
\def\AGammaL{{\sf A\Gamma L}}

 \def\PSU{{\sf PSU}}  
 \def\F{{\sf F}} \def\D{{\sf D}}

\def\GammaL{{\sf \Gamma L}}

\def\Sz{{\sf Sz}}

\newtheorem{theorem}{Theorem}[section]%
\newtheorem{lemma}[theorem]{Lemma}%
\newtheorem{example}[theorem]{Example}%

\input amssym.def
\input amssym.tex

\begin{document}

\title[Pentavalent graphs]{Pentavalent symmetric graphs of order four times an odd square-free integer}
\thanks{2000 MR Subject Classification 05C25, 05E18, 20B25.}
\thanks{This work was partially supported by the National Natural Science Foundation of China (11301468,11231008,11461004)}
\thanks{$^*$Corresponding author. E-mails: bolinggxu@163.com (B. Ling).}
\author[Ling, Lou, and Wu]{Bo Ling$^{a,*}$, Ben Gong Lou$^b$, and Ci Xuan Wu$^c$}
\address{a: School of Mathematics and Computer Sciences\\
Yunnan Minzu University\\
Kunming, Yunnan 650504, P.R. China}
\email{bolinggxu@163.com (B. Ling)}
\address{b: School of Mathematics and Statistics\\ Yunnan University,
Kunmin 650031, P. R. China}
\email{bengong188@163.com (B.G. Lou)}
\address{c: School of Statistics and Mathematics\\
Yunnan University of Finance and Economics \\
Kunming, Yunnan, P. R. China}
\email{wucixuan@gmail.com (C.X. Wu)}

\date\today

\begin{abstract}
A graph is said to be symmetric if its automorphism group is transitive on its arcs.
Guo et al. (Electronic J. Combin. 18, \#P233, 2011) and Pan et al. (Electronic J. Combin. 20, \#P36, 2013)
determined all pentavalent symmetric graphs of order $4pq$.
In this paper, we shall generalize this result by determining all connected pentavalent symmetric graphs of order four times an odd square-free integer.
It is shown in this paper that, for each of such graphs $\Ga$, either the full automorphism group $\Aut\Ga$ is isomorphic to $\PSL(2,p)$, $\PGL(2,p)$, $\PSL(2,p){\times}\ZZ_2$ or $\PGL(2,p){\times}\ZZ_2$, or $\Ga$ is isomorphic to one of 8 graphs.
\vskip0.1in


\end{abstract}
\maketitle

\section{Introduction}
All graphs in this paper are assumed to be finite, simple, connected and undirected.
Let $\Ga$ be a graph and denote $V\Ga$ and $A\Ga$ the vertex set and arc set of $\Ga$, respectively.
Let $G$ be a subgroup
of the full automorphism group $\Aut\Ga$ of $\Ga$.
Then $\Ga$ is called {\it $G$-vertex-transitive} and
{\it $G$-arc-transitive} if $G$ is transitive on $V\Ga$ and $A\Ga$, respectively.
An arc-transitive graph is also called {\it symmetric}.
It is well known that $\Ga$ is $G$-arc-transitive if and only if
$G$ is transitive on $V\Ga$ and the stabilizer $G_{\a}:=\{g\in G\mid \a^g=\a\}$
for some $\a\in V\Ga$ is transitive on the neighbor set $\Ga(\a)$ of $\a$ in $\Ga$.

The cubic and tetravalent graphs have been studied extensively in the literature.
In recent years, attention has moved on to pentavalent symmetric graphs
and a series of results have been obtained.
For example, all the possibilities of vertex stabilizers
of pentavalent symmetric graphs are determined in \cite{Guo-Feng,Zhou-Feng}.
Also, for distinct primes $p$, $q$ and $r$, the classifications of pentavalent symmetric graphs of order
$2pq$ and $2pqr$ are presented in \cite{Hua-1,YFD16}, respectively.
A classification of $1$-regular pentavalent graph
(that is, the full automorphism group acts regularly on its arc set)
of square-free order is presented in \cite{Li-Feng}.
Recently, pentavalent symmetric graphs of square-free order have been completely classified in \cite{LZW16}.
Furthermore, some classifications of pentavalent symmetric graphs of cube-free order also have been obtained in recent years.
For example, the classifications of pentavalent symmetric graphs of order $12p$, $4pq$ and $2p^2$ are presented
in \cite{Guo,PLL13,FZL16}.
The main purpose of this paper is to extend the results in \cite{Guo,PLL13} to four times an odd square-free integer case.

The main result of this paper is the following theorem.

\begin{theorem}\label{square-free order}
Let $n$ be an odd square-free integer and
let $\Ga$ be a pentavalent symmetric graph of order $4n$.
If $n$ has at least three prime factors,
then one of the following statements holds.
\begin{itemize}
\item[(1)] $\Aut\Ga\cong\PSL(2,p)$, $\PGL(2,p)$, $\PSL(2,p)\times\ZZ_2$ or $\PGL(2,p)\times\ZZ_2$, where $p\ge 29$ is a prime.
\item[(2)] The triple $(\Ga, n, \Aut\Ga)$ lies in the following Table \ref{table-4square-free}.
\begin{table}[ht]
 \[\begin{array}{|ccccccc|}
      \hline
      \ \Row& \Ga & n  &\Aut\Ga &(\Aut\Ga)_{\alpha} &\Transitivity&\Bipartite?\\
      \hline
    \ 1&\mathcal{C}_{17556}^{1} & 3{\cdot}7{\cdot}11{\cdot}19 &\J_1 &\D_{10}  &1-\transitive&\No \\
     \ 2&\mathcal{C}_{17556}^{2} & 3{\cdot}7{\cdot}11{\cdot}19 &\J_1 &\D_{10}  &1-\transitive&\No \\
      \ 3&\mathcal{C}_{17556}^{3} & 3{\cdot}7{\cdot}11{\cdot}19 &\J_1 &\D_{10}  &1-\transitive&\No \\
       \ 4&\mathcal{C}_{17556}^{4} & 3{\cdot}7{\cdot}11{\cdot}19 &\J_1 &\D_{10}  &1-\transitive&\No \\
       \ 5&\mathcal{C}_{17556}^{5} & 3{\cdot}7{\cdot}11{\cdot}19 &\J_1 &\D_{10}  &1-\transitive&\No \\
       \ 6&\mathcal{C}_{5852} & 7{\cdot}11{\cdot}19 &\J_1{\times}\ZZ_2 &\A_{5}  &2-\transitive&\Yes\\
           \ 7&\mathcal{C}_{780}^{1} & 3{\cdot}5{\cdot}13 &\PSL(2,25){\times}\ZZ_2 &\F_{20}  &2-\transitive&\No \\
                     \ 8&\mathcal{C}_{780}^{2} & 3{\cdot}5{\cdot}13 &\PSL(2,25){\times}\ZZ_2 &\F_{20}  &2-\transitive&\No \\
      \hline
\end{array}\]
\caption{\label{table-4square-free}}
\end{table}
\end{itemize}
\end{theorem}
\noindent{\bf Remark 1.1.} \begin{itemize}
\item[(a)] The graphs in Table \ref{table-4square-free} are introduced in Example \ref{exam-PSL(2,25)-J1}.

\item[(b)] It seems not feasible to determine all the possible values of $p$
in part (1) for general odd square-free integer $n$. However, if the number of the prime
divisors of $n$ is fixed, then it is not difficult to determine the possible values of $p$
and hence all corresponding graphs $\Ga$.
\end{itemize}
\section{Preliminaries}
We now give some necessary preliminary results.
The first one is a property of the Fitting subgroup,
see \cite[P. 30, Corollary]{Suzuki}.

\begin{lemma}\label{Fitting-sg}
Let $F$ be the Fitting subgroup of a group $G$. If $G$ is soluble, then $F\ne 1$ and
the centralizer $C_G(F)\le F$.
\end{lemma}

The maximal subgroups of $\PSL(2,q)$ are known, see \cite[Section 239]{Dickson}.

\begin{lemma}\label{Subg-PSL(2,q)}
Let $T=\PSL(2,q)$, where $q=p^n\ge 5$ with $p$ a prime.
Then a maximal subgroup of $T$ is isomorphic to one of the following groups,
where $d=(2,q-1)$.

\begin{enumerate}
\item[(1)] $\D_{2(q-1)/d}$, where $q\ne 5,7,9,11$;
\item[(2)] $\D_{2(q+1)/d}$, where $q\ne 7,9$;
\item[(3)] $\ZZ_{p}^n{:}\ZZ_{{(q-1)/d}}$;
\item[(4)] $\A_4$, where $q=p=5$ or $q=p\equiv 3,13,27,37~(\mod40)$;
\item[(5)] $\S_4$, where $q=p\equiv \pm 1~(\mod8)$
\item[(6)] $\A_5$, where $q=p\equiv\pm 1~(\mod5)$, or $q=p^2\equiv -1~(\mod5)$
with $p$ an odd prime;
\item[(7)] $\PSL(2,p^m)$ with $n/m$ an odd integer;
\item[(8)] $\PGL(2,p^{n/2})$ with $n$ an even integer.
\end{enumerate}
\end{lemma}

By \cite[Theorem 2]{COT}, we may easily derive the maximal subgroups of $\PGL(2, p)$.
\begin{lemma}\label{Subg-PGL(2,q)}
Let $T = \PGL(2, p)$ with $p\ge5$ a prime. Then a maximal subgroup of
$T$ is isomorphic to one of the following groups:
\begin{enumerate}
\item[(1)] $\ZZ_p:\ZZ_{p-1}$;
\item[(2)] $\D_{2(p+1)}$;
\item[(3)] $\D_{2(p-1)}$, where $p\ge7$;
\item[(4)] $\S_4$, where $p\equiv \pm 3~(\mod8)$;
\item[(5)] $\PSL(2,p)$.
\end{enumerate}
\end{lemma}

For a graph $\Ga$ and a positive integer $s$, an $s$-{\it arc} of $\Ga$ is
a sequence $\a_0,\a_1,\dots,\a_s$ of vertices such that
$\a_{i-1},\a_i$ are adjacent for $1\le i\le s$ and $\a_{i-1}\not=\a_{i+1}$
for $1\le i\le s-1$. In particular, a $1$-arc is just an arc. Then
$\Ga$ is called {\it $(G,s)$-arc-transitive} with $G\le\Aut\Ga$
if $G$ is transitive on the set of $s$-arcs of $\Ga$.
A $(G,s)$-arc-transitive graph is called {\it $(G,s)$-transitive}
if it is not $(G,s+1)$-arc-transitive.
In particular, a graph $\Ga$ is simply called {\it $s$-transitive}
if it is $(\Aut\Ga,s)$-transitive.

Let $\F_{20}$ denote the Frobenius group of order $20$. 
The following lemma determines the stabilizers of pentavalent symmetric graphs,
refer to \cite{Guo-Feng,Zhou-Feng}.

\begin{lemma}\label{weiss}
Let $\Ga$ be a pentavalent $(G,s)$-transitive graph, where $G\le\Aut\Ga$ and $s\ge 1$.
Let $\a\in V\Ga$. Then one of the following holds.
\begin{itemize}
\item[(a)] If $G_{\a}$ is soluble, then $s\le 3$ and $|G_{\a}|\div 80$.
Further, the couple $(s,G_{\a})$ lies in the following table.

\[\begin{array}{|c|c|c|c|} \hline
s & 1 & 2 & 3 \\ \hline
G_{\a} &  \ZZ_5,~\D_{10},~\D_{20}& \F_{20},~\F_{20}\times\ZZ_2 & \F_{20}\times\ZZ_4  \\ \hline
\end{array}\]
\vskip0.2in

\item[(b)] If $G_{\a}$ is insoluble, then $2\le s\le 5$, and $|G_{\a}|\div 2^{9}\cdot 3^2\cdot 5$.
Further, the couple $(s,G_{\a})$ lies in the following table.

\[\begin{array}{|c|c|c|c|c|} \hline
s & 2 & 3 & 4 & 5 \\ \hline
G_{\a} &  \A_5,\S_5 &  \A_4\times\A_5,(\A_4\times\A_5){:}\ZZ_2,& \ASL(2,4),\AGL(2,4),& \ZZ_2^6{:}\GammaL(2,4)  \\
       & &  \S_4\times\S_5 & \ASigmaL(2,4),\AGammaL(2,4)&  \\ \hline
       |G_{\a}|&60,120&720,1440,2880&960,1920,2880,5760&23040\\ \hline
\end{array}\]
\end{itemize}
\end{lemma}
\vskip0.2in

From \cite[pp. 134-136]{Gorenstein}, we can obtain the following lemma
by checking the orders of nonabelian simple groups.

\begin{lemma}\label{square-free-sg}
Let $n$ be an odd square-free integer such that $n$ has at least three prime factors.
Let $T$ be a nonabelian simple group of order $2^{i}\cdot3^{j}\cdot 5 \cdot n$,
where $1\leq i \leq 11$ and $0 \leq j \leq 2$. Let $p$ be the largest prime factor of $n$.
Then $T$ is listed Table \ref{table-simple-groups}.
\begin{table}[ht]
\[\begin{array}{|llllll|} \hline
T & |T| & n & T & |T| & n  \\ \hline
\M_{22} & 2^7{\cdot}3^2{\cdot}5{\cdot}7{\cdot}11 & 3{\cdot}7{\cdot}11 & \PSp(4,4) & 2^8{\cdot}3^2{\cdot}5^2{\cdot}17 & 3{\cdot}5{\cdot}17 \\
\M_{23}  & 2^7{\cdot}3^2{\cdot}5{\cdot}7{\cdot}11{\cdot}23 & 7{\cdot}11{\cdot}23& \PSL(2,25) & 2^3{\cdot}3{\cdot}5^2{\cdot}13 & 3{\cdot}5{\cdot}13  \\
\J_1 & 2^3{\cdot}3{\cdot}5{\cdot}7{\cdot}11{\cdot}19 & 7{\cdot}11{\cdot}19 & \PSL(2,2^8) & 2^8{\cdot}3{\cdot}5{\cdot}17{\cdot}257 & 3{\cdot}17{\cdot}257 \\
\J_2 & 2^7{\cdot}3^3{\cdot}5^2{\cdot}7 & 3{\cdot}5{\cdot}7 & \PSL(5,2) & 2^{10}{\cdot}3^2{\cdot}5{\cdot}7{\cdot}31 &  3{\cdot}7{\cdot}31 \\
\Sz(32) & 2^{10}{\cdot}5^2{\cdot}31{\cdot}41 & 5{\cdot}31{\cdot}41 &\PSL(2,2^6) & 2^6{\cdot}3^2{\cdot}5{\cdot}7{\cdot}13 & 3{\cdot}7{\cdot}13 \\
\PSU(3,4) & 2^6{\cdot}3{\cdot}5^2{\cdot}13 & 3{\cdot}5{\cdot}13 &&& \\ \hline
\M_{23}  & 2^7{\cdot}3^2{\cdot}5{\cdot}7{\cdot}11{\cdot}23 & 3{\cdot}7{\cdot}11{\cdot}23&
\J_1 & 2^3{\cdot}3{\cdot}5{\cdot}7{\cdot}11{\cdot}19 & 3{\cdot}7{\cdot}11{\cdot}19 \\
\M_{24}  & 2^{10}{\cdot}3^3{\cdot}5{\cdot}7{\cdot}11{\cdot}23 & 3{\cdot}7{\cdot}11{\cdot}23 &&&     \\ \hline
\PSL(2,p) & \frac{p(p+1)(p-1)}{2}\ (p\ge29) &&&& \\ \hline
\end{array}\]\caption{\label{table-simple-groups}}
\end{table}
\end{lemma}
\proof
If $T$ is a sporadic simple group, by \cite[P. 135-136]{Gorenstein}, $T =\M_{22}$, $\M_{23}$, $\M_{24}$, $\J_1$ or $\J_2$.
If $T=\A_n$ is an alternating group, since $3^4$ does not divide $|T|$,
we have $n\le8$, it then easily exclude that $T=\A_5$, $\A_6$, $\A_7$ or $\A_8$.
Hence no $T$ exists for this case.

Suppose now $T = X(q)$ is a simple group of Lie type,
where $X$ is one type of Lie groups, and $q ={r}^d$ is a prime power.
If $r\ge5$, as $|T|$ has at most three 3-factors, two 5-factors and
one $p$-factor, it easily follows from \cite[P. 135]{Gorenstein} that
the only possibility is $T = \PSL(2,p)$ with $p\ge29$ (note that $\PSL(2,p)$ with $5\le p \le23$
does not satisfy the condition of the lemma) or $\PSL(2,25)$, where $p$ is the largest prime factor of $n$.
If $r\le3$, as $2^{12}$ and $3^4$ do not divide $|T|$, then we have $T = \Sz(32)$, $\PSU(3,4)$, $\PSp(4,4)$, $\PSL(2,2^6)$, $\PSL(2,2^8)$ or $\PSL(5,2)$.
\qed

\vskip0.1in

A typical method for studying vertex-transitive graphs is taking normal quotients.
Let $\Ga$ be a $G$-vertex-transitive graph,
where $G\le\Aut\Ga$.
Suppose that $G$ has a normal subgroup $N$ which is intransitive on $V\Ga$.
Let $V\Ga_N$ be the set of $N$-orbits on $V\Ga$.
The {\it normal quotient graph} $\Ga_N$ of $\Ga$ induced by $N$ is defined as the graph
with vertex set $V\Ga_N$, and $B$ is adjacent to $C$ in $\Ga_N$
if and only if there exist vertices $\b\in B$ and $\g\in C$ such that
$\b$ is adjacent to $\g$ in $\Ga$.
In particular, if $\val(\Ga)=\val(\Ga_N)$, then $\Ga$ is called a {\it normal cover } of $\Ga_N$.

A graph $\Ga$ is called {\it $G$-locally primitive} if,
for each $\a\in V\Ga$, the stabilizer $G_{\a}$ acts primitively on $\Ga(\a)$.
Obviously, a pentavalent symmetric graph is locally primitive.
The following theorem gives a basic method for studying vertex-transitive locally primitive graphs,
see \cite[Theorem 4.1]{Praeger} and \cite[Lemma 2.5]{L-Pan}.

\begin{theorem}\label{praeger}
Let $\Ga$ be a $G$-vertex-transitive locally primitive graph, where $G\le\Aut\Ga$,
and let $N\lhd G$ have at least three orbits on $V\Ga$. Then the following statements hold.
\begin{itemize}
\item[(i)] $N$ is semi-regular on $V\Ga$, $G/N\le\Aut\Ga_N$,
and $\Ga$ is a normal cover of $\Ga_N$;
\item[(ii)] $G_{\a}\cong (G/N)_{\g}$, where $\a\in V\Ga$ and $\g\in V\Ga_N$;
\item[(iii)] $\Ga$ is $(G,s)$-transitive if and only if $\Ga_N$
is $(G/N,s)$-transitive, where $1\le s\le 5$ or $s=7$.
\end{itemize}
\end{theorem}
\vskip0.1in


For reduction, we need some information of pentavalent symmetric graphs of order $4pq$,
stated in the following lemma, see \cite[Theorem 4.1]{Guo} and \cite[Theorem 3.1]{PLL13}.

\begin{lemma}\label{4pq-gps}
Let $\Ga$ be a pentavalent symmetric graph of order $4pq$, where $q> p\ge 3$ are primes.
Then the couple $(\Aut\Ga,(\Aut\Ga)_{\a})$ lies in the following Table $3$, where $\a\in V\Ga$.


\begin{table}[ht]
\[\begin{array}{|llll|} \hline
\Ga & (p,q) &  \Aut\Ga & (\Aut\Ga)_{\a}  \\ \hline
 \mathcal{C}_{60}&(3,5)&\A_5{\times}\D_{10}&\D_{10}\\
 \mathcal{C}_{132}^{1}&(3,11)&\PSL(2,11){\times}\ZZ_2&\D_{10}\\
 \mathcal{C}_{132}^{i}, 2\le i\le 4&(3,11)&\PGL(2,11)&\D_{10}\\
 \mathcal{C}_{132}^{5}&(3,11)&\PGL(2,11){\times}\ZZ_2&\D_{20}\\
 \CC^{(2)}_{574} & (7,41) & \PSL(2,41)\times\ZZ_2 & \A_5  \\
 \CC_{4108} & (13,79) & \PSL(2,79) & \A_5 \\ \hline
\end{array}\]
\caption{\label{table-4pq}}
\end{table}
\end{lemma}
\noindent{\bf Remark 2.8.}
\begin{itemize}
\item[(a)] Suppose that $\Ga$ is one of the graphs in Lemma \ref{4pq-gps} and $M$ is an arc-transitive subgroup of $\Aut\Ga$.
 Then $M$ is insoluble (for convenience, we prove this conclusion in Lemma \ref{Lemma-1} and we remark
that Lemma \ref{Lemma-1} is independent where it is used).
\item[(b)] By Magma \cite{Magma}, the graphs $\mathcal{C}_{66}^{(2)}$ and $\mathcal{C}_{132}^{5}$ in \cite[Theorem 4.1]{Guo} are isomorphic,
$\Aut(\mathcal{C}_{132}^{5})\cong\PGL(2,11){\times}\ZZ_2$.
\end{itemize}
\vskip0.1in

The final lemma of this section gives some information about the pentavalent symmetric graphs of square-free order,
refer to \cite[Theorem 1.1]{YFD16} and \cite[Theorem 1.1]{LZW16}.
\begin{lemma}\label{square-free-gras}
Let $\Ga$ be a pentavalent symmetric graph of order $2n$,
where $n$ is an odd square-free integer and has at least three prime factors. Then one of
the following statements holds.
\begin{itemize}
\item[(1)] $\Aut\Ga$ is soluble and $\Aut\Ga\cong\D_{2n}:\ZZ_5$.
\item[(2)] $\Aut\Ga=\PSL(2,p)$ or $\PGL(2,p)$, where $p\ge 5$ is a prime.
\item[(3)] The triple $(\Ga, 2n, \Aut\Ga)$ lies in the following Table \ref{table-square-free}.
\begin{table}[ht]
 \[\begin{array}{|llll|}
     \hline
      \ \Ga & 2n  &\Aut\Ga &(\Aut\Ga)_{\alpha} \\
      \hline
   \ \mathcal{C}_{390} & 390 &\PSL(2,25) &\F_{20} \\
   \ \mathcal{C}_{2926} & 2926 &\J_1 &\A_5  \\
      \hline
\end{array}\]
\caption{\label{table-square-free}}
\end{table}
\end{itemize}
\end{lemma}

\section{Some examples}

In this section, we give some examples of pentavalent symmetric graphs
of order $4n$ with $n$ an odd square-free integer.

For a given small permutation group $X$, we may determine all graphs which admit $X$
as an arc-transitive automorphism group by using Magma \cite{Magma}.
It is then easy to have the following result.
\begin{example}\label{exam-PSL(2,25)-J1}\rm
{\rm(1)} There is a unique pentavalent symmetric graph of order $5852$
which admits $\J_1\times\ZZ_2$ as an arc-transitive automorphism group;
and its full automorphism group is $\J_1\times\ZZ_2$.
This graph is denoted by $\CC_{5832}$ which satisfies the conditions in Row 6 of Table \ref{table-4square-free}.

{\rm(2)} There are five pentavalent symmetric graphs of order $17556$
admitting $\J_1$ as an arc-transitive automorphism group;
and their full automorphism group are all isomorphic to $\J_1$.
These five graphs are denoted by $\CC_{17556}^i$ which satisfy the conditions in Row 1 to Row 5 of Table \ref{table-4square-free}, where $1\le i\le5$.

{\rm(3)} There are two pentavalent symmetric graphs of order $780$
which admits $\PSL(2,25)\times\ZZ_2$ as an arc-transitive automorphism group;
and their full automorphism group are all isomorphic to $\PSL(2,25)\times\ZZ_2$.
These two graphs are denoted by $\CC_{780}^j$ which satisfy the conditions in Row 7 to Row 8 of Table \ref{table-4square-free}, where $1\le j\le2$.
\end{example}

\section{Proof of Theorem \ref{square-free order}}
Let $n$ be an odd square-free integer and $n$ has at least three prime factors.
Let $\Ga$ be a pentavalent symmetric graph of order $4n$.
Set $\A=\Aut\Ga$. By Lemma~\ref{weiss}, $|\A_{\a}|\div 2^{9}\cdot 3^{2}\cdot 5$,
and hence $|\A|\div 2^{11}\cdot 3^{2}\cdot 5\cdot n$.
Assume that $n=p_1 p_2\cdots p_s$,
where $s\ge3$ and ${p_i}'s$ are distinct primes.

We first consider the case where $\A$ is soluble.

\begin{lemma}\label{soluble}
 Assume that $\A$ is soluble. Then no graph $\Ga$ exists.
 \end{lemma}
\proof Let $F$ be the Fitting subgroup of $\A$.
By Lemma~\ref{Fitting-sg},
$F\ne 1$ and $\C_{\A}(F)\le F$.
Further, $F=\O_2(\A)\times\O_{p_1}(\A)\times\O_{p_2}(\A)\times\cdots\times\O_{p_s}(\A)$,
where $\O_2(\A)$, $\O_{p_1}(\A)$, $\O_{p_2}(\A)$,$\dots$,$\O_{p_s}(\A)$
denote the largest normal $2$-, $p_1$-, $p_2$-, $\dots$, $p_s$-subgroups of $\A$, respectively.

For each $p_i\in\{p_1,p_2,\dots,p_s\}$, $\O_{p_i}(\A)$ has at least three orbits on $V\Ga$,
by Theorem \ref{praeger}, $\O_{p_i}(\A)$ is semi-regular on $V\Ga$.
Therefore, $F$ is semi-regular on $V\Ga$ and $\O_{p_i}(\A)\le \ZZ_{p_i}$.
This argument also proves $\O_{2}(\A)\le\ZZ_{4}$ or $\ZZ_2^2$. If $\O_2(\A)=\ZZ_4$ or $\ZZ_2^2$,
then by Theorem \ref{praeger}, the normal quotient graph $\Ga_{\O_2(\A)}$ is a pentavalent symmetric graph of
odd order, which is a contradiction.
Thus, $\O_2(\A)\le\ZZ_2$, $F\cong\ZZ_m$, where $m\div 2n$.
It implies that $\C_{\A}(F)\ge F$, and so $\C_{\A}(F)=F$.

If $F$ has at least three orbits on $V\Ga$, then, by Theorem \ref{praeger},
$\Ga_F$ is $\A/F$-arc-transitive. Since $\A/F=\A/\C_{\A}(F)\le\Aut(F)$ is abelian,
we have $(\A/F)_{\delta}=1$, where $\delta\in V\Ga_F$, which is a contradiction.

Thus, $F$ has at most two orbits on $V\Ga$.
If $F$ is transitive on $V\Ga$,
then $F$ is regular on $V\Ga$.
It follows that $\Ga\cong\Cay(F,S)$ is a normal arc-transitive Cayley graph of $F$.
Then we easily conclude that $S$ consists of the involutions of $F$. Since $F$ has at most one involution,
it follows that $F=\l S \r\le\ZZ_2$,
which is a contradiction.

Hence $F$ has two orbits on $V\Ga$ and $F\cong\ZZ_{2n}$ and
$K := \O_{p_3}(\A)\times\O_{p_4}(\A)\times\O_{p_s}(\A)\cong\ZZ_{p_3p_4\ldots p_s}$.
Since $K\unlhd\A$ has $4p_1p_2$ orbits on $V\Ga$, by Theorem \ref{praeger}(i),
$\Ga_K$ is an $\A/K$-arc-transitive pentavalent graph of order $4p_1p_2$, and
hence $\Ga_K$ satisfies Lemma \ref{4pq-gps}. Since $\A/K$ is soluble, by Remark 2.8, a contradiction occurs. \qed

We next consider the case where $\A$ is insoluble.

\begin{lemma}\label{insoluble-1}
 Assume that $\A$ is insoluble and has no nontrivial soluble normal subgroup.
Then $\Aut\Ga\cong\J_1$, $\PSL(2,p)$ or $\PGL(2,p)$. Further, If $\Aut\Ga\cong\J_1$,
then $\Ga\cong \CC_{17556}^i$ which satisfy the conditions in
Row 2 to Row 5 of Table \ref{table-4square-free} of Theorem \ref{square-free order}, where $1\le i\le5$.
 \end{lemma}
\proof
Let $N$ be a minimal normal subgroup of $\A$.
Then $N=S^d$, where $S$ is a nonabelian simple group and $d\geq 1$.

If $N$ has more than three orbits on $V\Ga$, then by Theorem \ref{praeger}, $N$ is semi-regular on $V\Ga$ and so $|N|$ divides $4n$.
Since $N$ is a direct product of nonabelian groups, it implies that $4$ divides $|N|$.
Again by Theorem \ref{praeger}, $\Ga_N$ is a pentavalent symmetric graph of odd order $n$, a contradiction.
Hence, $N$ has at most two orbits on $V\Ga$,
so $2n$ divides $|N|$.

Moreover, since $p_s>5$, $p_s$ divides $|N|$ and $p_s^2$ does not divide $|N|$ as
$|\A|\div 2^{11}\cdot 3^{2}\cdot 5\cdot p_1p_2\cdots p_s$,
we conclude that $d=1$ and $N=S$ is a nonabelian simple group.
Let $C=\C_\A(S)$.
Then $C\lhd \A$, $C\cap S=1$ and $\l C,S\r=C\times S$.
If $C\not=1$, then $C$ is insoluble as $C\lhd \A$ and $\A$ has no soluble normal subgroup.
It follows that $4$ divides $|C|$. A similarly argument with the above paragraph, we have $2n$ divides $|C|$.
Hence $4n^2$ divides $|\A|=2^{11}\cdot 3^{2}\cdot 5\cdot n$,
and so $n$ divides $2^{9}\cdot 3^{2}\cdot 5$. It implies that $n=3\cdot5$,
a contradiction with $n$ having at least three prime factors.
So $C=1$ and $\A=\A/C\le\Aut(S)$, that is,
$\A$ is almost simple with socle $S$.

If $S_\a=1$, then $S$ acts regularly on $V\Ga$. Hence $S$ is a non-abelian simple group such that $|S|=4n$.
By checking the orders of nonabelian simple groups (see \cite[P. 135-136]{Gorenstein} for example),
we have that $S=\PSL(2,p)$ and so $\A\le\Aut(S)=\PGL(2,p)$, which is impossible as $\A$ is transitive on $A\Ga$, $|\A|\le 2|S|$ and $|A\Ga|=5|S|$.
Hence $S_\a\not=1$. Since $\Ga$ is connected and $S\lhd \A$, we have
$1\ne S_{\a}^{\Ga(\a)}\lhd \A_{\a}^{\Ga(\a)}$, it follows that $5\div |S_{\a}|$,
we thus have $10\cdot p_1p_2\cdots p_s$ divides $|S|$.

Thus, $\soc(\A)=S$ is a nonabelian simple group such that $|S|\div2^{11}\cdot3^2\cdot5\cdot n$ and $10\cdot n\div|S|$.
Hence the triple $(S, |S|, n)$ lies in Table 2 of Lemma \ref{square-free-sg}.
We will analyse all the candidates one by one in the following.

If $\cong\PSL(2, p)$ with $p\ge29$ a prime, then $\A\cong\PSL(2, p)$ or $\PGL(2, p)$,
the Lemma holds. If $(S, n) = (\J_1, 3{\cdot}7{\cdot}11{\cdot}19)$,
then $|V\Ga|=17556$ and $\A\cong \J_1$ as $\Out(\J_1) = 1$.
It then follows from Example \ref{exam-PSL(2,25)-J1} that
$\Ga\cong \CC_{17556}^i$ which satisfy the conditions in
Row 1 to Row 5 of Table \ref{table-4square-free} of Theorem \ref{square-free order}, where $1\le i\le5$.

Assume $(S, n) = (\Sz(32), 5{\cdot}31{\cdot}41)$. Since $\Out(\Sz(32))\cong\ZZ_5$ (see Atlas \cite{Atlas} for example),
$\A\cong\Sz(32)$ or $\Sz(32).\ZZ_5$, so $|\A_{\a}| = \frac{|\A|}{4n}= 1280$ or 6400, which is not possible by Lemma \ref{weiss}.
Similarly, for the case $(S, n) = (\PSL(5, 2), 3{\cdot}7{\cdot}31)$, then $\A\cong\PSL(5, 2)$ or $\PSL(5, 2).\ZZ_2$
as $\Out(\PSL(5, 2))\cong\ZZ_2$. Thus, $|\A_{\a}| =\frac{|\A|}{4n}=3840$ or 7680, which
is impossible by Lemma \ref{weiss}. For the case where $(S, n) = (\PSL(2, 2^8), 3{\cdot}17{\cdot}257)$,
since $\A\cong\PSL(2, 2^8).O$, where $O\le \Out(\PSL(2, 2^8))\cong\ZZ_8$, we have $|\A_{\a}| =\frac{|\A|}{4n}=2^k{\cdot}5$,
where $6\le k \le 9$, which is also impossible by Lemma \ref{weiss}. For the case where $(S,n)=(\PSU(3,4),3{\cdot}5{\cdot}13)$,
since $\A\cong\PSU(3,4).O$, where $O\le \Out(\PSU(3,4))\cong\ZZ_4$, we have $|\A_{\a}| =\frac{|\A|}{4n}=2^k{\cdot}5$,
where $4\le k \le 6$, which is impossible by Lemma \ref{weiss}.

Assume $(S,n)=(\PSp(4,4),3{\cdot}5{\cdot}17)$. Since $S\le\A\le\Aut(S)\cong\PSp(4, 4).\ZZ_4$,
we have $|\A_{\a}| =\frac{|\A|}{4n}=960$, 1920 or 3840. If $|\A_{\a}|=960$ or $1920$, then by Lemma \ref{weiss}, $\A_{\a}\cong\ASL(2,4)$ or $\ASigmaL(2,4)$.
However, by Atlas \cite{Atlas}, $\PSp(4,4)$ has no subgroup isomorphic to $\ASL(2,4)$ and $\PSp(4,4).\ZZ_2$ has no subgroup isomorphic to
$\ASigmaL(2,4)$. If $|\A_{\a}|=3840$, then also by Lemma \ref{weiss}, a contradiction occurs.

Assume $(S, n) = (\PSL(2, 2^6), 3{\cdot}7{\cdot}13)$. Recall that $S$ has at most two orbits on
$V\Ga$, $|S_{\a}| = \frac{|S|}{4n}=240$ or $\frac{|S|}{2n}=480$. However, by Lemma \ref{Subg-PSL(2,q)},
$\PSL(2, 2^6)$ has no maximal subgroup with order a multiple of 240, a contradiction occurs.
Similarly, for the case $(S,n)=(\J_2, 3{\cdot}5{\cdot}7)$. Then $|S_{\a}| = \frac{|S|}{4n}=2880$ or $\frac{|S|}{2n}=5760$.
By Atlas \cite{Atlas}, $\J_2$ has no maximal subgroup with order a multiple of 2880, a contradiction also occurs.

Assume $S\cong\M_{23}$. Then $n =3{\cdot}7{\cdot}11{\cdot}23$ or $7{\cdot}11{\cdot}23$, and as $\Out(\M_{23}) = 1$,
we have $\A = S$ and $|\A_{\a}| = \frac{|\M_{23}|}{4n}=480$ or $1440$. By Lemma \ref{weiss}, it is impossible for the case $|\A_{\a}|=480$.
For the latter case, by a direct computation using
Magma \cite{Magma}, no graph $\Ga$ exists. If $(S, n) = (\M_{22}, 7{\cdot}11{\cdot}23)$,
as $\Out(\M_{22})\cong\ZZ_2$, we have $\A\cong\M_{22}$ or $\M_{22}.\ZZ_2$, so $|\A_{\a}| = \frac{|\A|}{4n}=480$ or 960, a
computation by Magma \cite{Magma} shows that no graph $\Ga$ exists.
Similarly, we can exclude the case where $(S,n)=(\PSL(2,25),3{\cdot}5{\cdot}13)$ by Magma \cite{Magma}.

Finally, assume $(S, n) = (\M_{24}, 3{\cdot}7{\cdot}11{\cdot}23)$ or $(\J_1, 3{\cdot}7{\cdot}11{\cdot}19)$. Since $\Out(\M_{24})=\Out(\J_1) = 1$,
we always have $\A = S$. Hence $|\A_{\a}| = \frac{|\A|}{4n}= 11520$ or 10. A
computation by Magma \cite{Magma} also shows that no graph $\Ga$ exists.
\qed

We next assume that $\A$ has a nontrivial soluble normal subgroup.
Let $N$ be a minimal soluble normal subgroup of $\A$.
Then there exists a prime $r\div 4n$ such that $N\cong\ZZ_r^d$.
Further, $N$ has at least three orbits on $V\Ga$.
It follows from Theorem \ref{praeger} that $N$ is semi-regular on $V\Ga$, and so $|N|=|\ZZ_r|^d\div|V\Ga|=4n$.
If $d\ge2$, then $(r,d)=(2,2)$. It follows that $\Ga_N$ is an arc-transitive graph of odd order, a contradiction.
Hence $d=1$, $N=\ZZ_r$.
The next lemma consider the case where $r=2$.

\begin{lemma}\label{insoluble-2}
Assume that $\A$ is insoluble and has a minimal soluble normal subgroup $N=\ZZ_2$.
Then one of the following statements holds:
\begin{itemize}
\item[(1)] $\Aut\Ga\cong\PSL(2,p)\times\ZZ_2$ or $\PGL(2,p)\times\ZZ_2$, where $p\ge 29$ is a prime.
\item[(2)] $\Aut\Ga\cong\PSL(2,25)\times\ZZ_2$ and $\Ga$ is isomorphic to $\CC_{780}^i$ in Table 1, where $1\le i\le2$.
\item[(3)] $\Aut\Ga\cong\J_1\times\ZZ_2$ and $\Ga$ is isomorphic to $\CC_{5852}$ in Table 1.
\end{itemize}
 \end{lemma}
\proof
Since $N$ has more than three orbits on $V\Ga$, then by Theorem \ref{praeger},
$\Ga_N$ is an $\A/N$-arc-transitive pentavalent graph of order $\bar n=2n$.
It follows that $\Ga_N$ is isomorphic to one of the graphs in Lemma \ref{square-free-gras}.
Since $\A/N\le\Aut\Ga_N$ and $\A/N$ is insoluble, we have that $\Aut\Ga_N$ is insoluble and so
$\Aut\Ga_N\cong\PSL(2,p)$, $\PGL(2,p)$, $\PSL(2,25)$ or $\J_1$.
Let $\bar\A:=\Aut\bar\Ga$.

Suppose that $\bar\A\cong\PSL(2, p)$ or $\PGL(2, p)$.
Since $\A/N$ is insoluble, by Lemma \ref{Subg-PSL(2,q)} and Lemma \ref{Subg-PGL(2,q)},
$\A/N$ is isomorphic to $\A_5$, $\PSL(2,p)$ or $\PGL(2,p)$. Since $\A/N$ is transitive on $A\Ga_N$,
we can further conclude that $\A/N$ is isomorphic to $\PSL(2,p)$ or $\PGL(2,p)$.
Therefore, $\A\cong N.\PSL(2,p)$ or $N.\PGL(2,p)$, that is, $\A\cong\PSL(2,p){\times}\ZZ_2$,
$\SL(2,p)$, $\PGL(2,p){\times}\ZZ_2$ or $\SL(2,p).\ZZ_2$. Assume first that $\A\cong\SL(2,p)$.
Note that $\SL(2,p)$ has a unique central involution. Then by Lemma \ref{weiss}, $\A_{\a}\cong\ZZ_5$.
It follows that $|V\Ga| =|\A :\A_\a|$ is divisible by 8 as $|\SL(2,p)|$ is divisible by 8, a contradiction.
Assume next that $\A\cong\SL(2,p).\ZZ_2$. Then $\A$ contains a normal subgroup $H$ isomorphic to $\SL(2,p)$.
Since $8\div|H|$, we have $H_\a\not=1$. By Theorem \ref{praeger},
$H$ has at most two orbits on $V\Ga$ and so $\frac{|\A_{\a}|}{|H_{\a}|}\div2$.
If $H$ is transitive on $V\Ga$, then $H$ is arc-transitive. A similar argument with
the case $\A\cong\SL(2,p)$, a contradiction occurs. Therefore, $H$ has two orbits on $V\Ga$ and so
$H_{\a}=\A_{\a}$. Since $H$ has a unique central involution, by Lemma \ref{weiss}, $\A_{\a}\cong\ZZ_5$,
it follows that $|V\Ga|=|\A:\A_\a|$ is divisible by 16, a contradiction.
Therefore, $\A\cong\PSL(2,p){\times}\ZZ_2$ or $\PGL(2,p){\times}\ZZ_2$ in this case.

Suppose that $\bar\A\cong\PSL(2,25)$. Since $\Ga_N$ is $\A/N$-arc-transitive, we have that $5{\cdot}390\div|\A/N|$.
By checking the maximal subgroup of $\PSL(2,25)$ (see Atlas \cite{Atlas} for example), we have that $\A/N=\bar\A\cong\PSL(2,25)$.
It follows that $\A\cong\SL(2,25)$ or $\PSL(2,25){\times}\ZZ_2$. If $\A\cong\PSL(2,25){\times}\ZZ_2$, then by Example \ref{exam-PSL(2,25)-J1},
$\Ga\cong\CC_{780}^i$ in Table 1, where $1\le i\le2$. If $\A\cong\SL(2,25)$, then by Magma \cite{Magma}, no graph exists.

Suppose that $\bar\A\cong\J_1$. Similarly, since $\Ga_N$ is $\A/N$-arc-transitive, we have that $5{\cdot}2926\div|\A/N|$.
By checking the maximal subgroup of $\J_1$ (see Atlas \cite{Atlas} for example), we have that $\A/N=\bar\A\cong\J_1$.
Since the Schur multiplier of $\J_1$ is $\ZZ_1$, $\A\cong N.\J_1\cong\J_1{\times}\ZZ_2$. By Example \ref{exam-PSL(2,25)-J1},
$\Ga\cong\CC_{5852}$ in Table 1. \qed

Finally, suppose that $r>2$. We first prove the following lemma.

\begin{lemma}\label{Lemma-1}
Assume that $\Sigma$ is isomorphic to one of
the graphs listed in Lemma \ref{4pq-gps}, in Lemma \ref{insoluble-1} and in Lemma \ref{insoluble-2}.
If $M$ is an arc-transitive subgroup of $\Aut\Sigma$, then $M$ contains a subgroup isomorphic to $\PSL(2,p)$, $\J_1$, $\PSL(2,25)$ or $\A_5$.
\end{lemma}
\proof By checking the graphs in Lemma \ref{4pq-gps}, in Lemma \ref{insoluble-1} and in Lemma \ref{insoluble-2},
we have that $\Aut\Sigma$ is isomorphic to one of the groups $\PSL(2,p)$, $\PGL(2,p)$, $\PSL(2,p){\times}\ZZ_2$,
$\PGL(2,p){\times}\ZZ_2$, $\J_1$, $\J_1{\times}\ZZ_2$, $\PSL(2,25){\times}\ZZ_2$ or $\A_5{\times}\D_{10}$.
If $\Aut\Sigma\cong\PSL(2,p)$, then since $p\div n$ and $20n\div|M|$, by Lemma \ref{Subg-PSL(2,q)},
$M\le\ZZ_p:\ZZ_{\frac{p-1}{2}}$ or $M=\Aut\Sigma\cong\PSL(2,p)$. If $M\le\ZZ_p:\ZZ_{\frac{p-1}{2}}$,
then $M\cong \ZZ_p:\ZZ_l$ for some $l\div \frac{p-1}{2}$. Thus, $M$ has a normal subgroup, say $S\cong\ZZ_p$,
which has more than three orbits on $V\Sigma$. It then follows from Theorem \ref{praeger} that the normal quotient graph
$\Sigma_S$ is $M/S$-arc-transitive, a contradiction occurs as $M/S\cong\ZZ_l$ is cyclic.
Hence, $M\not\le \ZZ_p:\ZZ_{\frac{p-1}{2}}$ and so $M=\Aut\Sigma\cong\PSL(2,p)$.
If $\Aut\Sigma\cong\PGL(2,p)$, then since $20n\div|M|$, by Lemma \ref{Subg-PGL(2,q)},
$\M\le\ZZ_p:\ZZ_{p-1}$, $M\le\PSL(2,p)$ or $M=\Aut\Sigma\cong\PGL(2,p)$.
A similar argument, we can conclude that $M\ge\PSL(2,p)$. Similarly,
we can further show that $M\ge\PSL(2,p)$ for the case $\Aut\Sigma\cong\PSL(2,p){\times}\ZZ_2$ or $\PGL(2,p){\times}\ZZ_2$.

If $\Aut\Sigma\cong\J_1$, then since $\Sigma$ is isomorphic to $\CC_{2926}$, we have that $5{\cdot}2926\div|M|$.
By checking the maximal subgroup of $\J_1$ (see Atlas \cite{Atlas}), we have that $M=\Aut\Sigma\cong\J_1$.
We can further show that $M\ge\J_1$ for the case $\Aut\Sigma\cong\J_1{\times}\ZZ_2$,
$M\ge\PSL(2,25)$ for the case $\Aut\Sigma\cong\PSL(2,25){\times}\ZZ_2$ and $M\ge\A_5$ for the case $\Aut\Sigma\cong\A_5{\times}\D_{10}$. \qed

Now assume that $\A$ has a minimal soluble normal subgroup $N=\ZZ_r$ for $r>2$.
\begin{lemma}\label{Lemma2}
Assume that $\A$ has a minimal soluble normal subgroup $N=\ZZ_r$ for $r>2$.
Then the normal quotient $\Ga_N$ is not isomorphic to any of the graphs $\Sigma$ listed in Lemma \ref{Lemma-1}.
\end{lemma}
\proof Suppose to the contrary that $\Ga_N$ is isomorphic to one of the graphs in Lemma \ref{Lemma-1}.
Then by Theorem \ref{praeger}, $\A/N\le\Aut\Ga_N$ is transitive on $A\Ga_N$.
Let $\Omega:=\{\PSL(2,p), \J_1, \PSL(2,25), \A_5\}$.
It follows from Lemma \ref{Lemma-1} that there exists a subgroup $M/N$ of $\A/N$ isomorphic to one of the groups in $\Omega$.
Since now $M/N\le\A/N\le\Aut\Ga_N$, it follows from the structure of $\Aut\Ga_N$ that $M/N\unlhd\A/N$.
Therefore, $M'{\sf char} M\unlhd\A$, it implies that $M'\unlhd\A$.
On the other hand, since the order of the Schur multiplier of a group in
$\Omega$ is less than or equal to 2 (see \cite[Theorem 7.1.1]{Kar87} for $\PSL(2,p)$ and Atlas \cite{Atlas} for the others) and $r>2$,
we have that $M'\in\Omega$ and $4\div|M'|$.
If $M'$ has more than three orbits on $V\Ga$, then by Theorem \ref{praeger},
$\Ga_{M'}$ is a pentavalent symmetric graph of odd order, a contradiction.
Thus, $M'$ has at most two orbits on $V\Ga$ and so $2n$ divides $|M'|$. Let $\bar\A:=\Aut\Ga_N$, $\bar n:=\frac{n}{r}$ and $\bar M:=M/N$.
Then $M'\cong\bar M$.

Assume first that $r>5$. Then since $|\bar M|\div|\bar\A|\div2^{11}{\cdot}3^2{\cdot}5{\cdot}\frac{n}{r}$
and $n$ is an odd square-free integer, we have that $r$ does not divide $|\bar M|=|M'|$.
It implies that $M'$ has at least $r$ orbits on $V\Ga$, a contradiction.

Assume next that $r=3$.
Since $\bar M\cong M'$ has at most two orbits on $V\Ga_N$ (if not $(\Ga_N)_{\bar M}$
is a pentavalent symmetric graph of odd order, a contradiction),
we have that $|\bar M:\bar M_{\delta}|=2\bar n$ or $4\bar n$, where $\delta\in V\bar\Ga$.
Now $2n$ divides $|\bar M|$ and $|\bar M:\bar M_{\delta}|=\frac{2n}{r}$ or $\frac{4n}{r}$.
It implies that $r=3$ divides $\bar M_{\delta}$.
Therefore $3\div|\bar\A_{\delta}|$.
By Lemma \ref{weiss}, $\bar\A_{\delta}$ is nonsolvable,
because $|\bar\A_{\delta}|$ does not divide $80$,
forcing that $\bar M_{\delta}$ is nonsolvable. If $\bar M\cong\PSL(2,p)$,
then by Lemma \ref{Subg-PSL(2,q)},
$\bar M_{\delta}\cong\A_5$.
Hence $M'_{\a}\le(M'N)_{\a}\cong(M'N/N)_{\delta}=\bar M_{\delta}\cong\A_5$ by Theorem \ref{praeger} (ii).
Note that $|M'_{\a}|=20$, it contradicts that $\A_5$ has no subgroup of order $20$.
If $\bar M\cong\J_1$, then $\Ga_N\cong\CC_{5852}$ or $\CC_{17556}^i$ in Table \ref{table-4square-free}, where $1\le i\le5$.
If $\Ga_N\cong\CC_{17556}^i$, then $\bar\A_{\delta}\cong\D_{10}$ is soluble, a contradiction.
If $\Ga_N\cong\CC_{5852}$, then $\bar M_{\delta}=\bar\A_{\delta}\cong\A_5$.
A similar argument with the case $\bar M\cong\PSL(2,p)$ leads to a contradiction.
If $\bar M\cong\A_5$, then $\Ga_N\cong\CC_{60}$ in Table 3 and $\bar\A_{\delta}\cong\D_{10}$ is soluble, a contradiction.
If $\bar M\cong\PSL(2,25)$, then $\Ga_N\cong\CC_{780}^1$ or $\CC_{780}^2$ in Table \ref{table-4square-free}
and $\bar\A_{\delta}\cong\F_{20}$ is soluble, also a contradiction.

Finally assume that $r=5$.
Then $|M':M'_{\a}|=2n$ or $4n$
as $M'$ has at most two orbits on $V\Ga$.
Since $\Ga$ is connected and $1\ne M'_{\a}\lhd \A_{\a}$, we have
$1\ne {M'_{\a}}^{\Ga(\a)}\lhd \A_{\a}^{\Ga(\a)}$, it follows that $5\div|M'_{\a}|$.
On the other hand, since $\bar M$ has at most two orbits on $V\Ga_N$,
we have that $|\bar M:\bar M_{\delta}|=2\bar n$ or $4\bar n$.
Note that $\bar M\cong M'$ and $r=5$.
Hence $\frac{|\bar M_{\delta}|}{|M'_{\a}|}=5$,
it follows that $5^2\div|\bar M_{\delta}|\div|\bar\A_{\delta}|$,
a contradiction with $|\bar\A_{\delta}|\div2^9\cdot3^2\cdot5$ by Lemma \ref{weiss}, a contradiction.
\qed
The final lemma completes the proof of Theorem \ref{square-free order}.
\begin{lemma}\label{Lemma3}
Assume $\A$ is insoluble. Then $\A$ has no minimal soluble normal
subgroup isomorphic to $\ZZ_r$ with $r>2$.
\end{lemma}
\proof
Suppose that, on the contrary, $\A$ has a minimal soluble normal subgroup $N=\ZZ_r$ with $r>2$.
We prove the lemma by induction on the order of $\Ga$.

Assume first that $n=pqt$ has three prime factors (Note that, by Table 3, the conclusion of Lemma \ref{Lemma3} does not hold for $n=pq$).
Without loss of generality, we may assume that $r=t$. Then $\Ga_N$ is a pentavalent symmetric graph of order $4pq$.
By Lemma \ref{4pq-gps}, $\Ga_N$ is isomorphic to one of the graphs in Table 3, which contradicts to Lemma \ref{Lemma2}.

Assume next that $n$ has at least four prime factors. Note that $\Aut\Ga_N$ is insoluble.
If $\Aut\Ga_N$ has no nontrivial soluble normal subgroup, then $\Ga_N$ is isomorphic to one of the graphs in Lemma \ref{insoluble-1},
which contradicts to Lemma \ref{Lemma2}. If $\Aut\Ga_N$ has a minimal soluble normal subgroup $\bar N$,
then we can also conclude that $\bar N\cong\ZZ_f$ with $f$ a prime. If $f>2$, then by induction, no such $\Ga_N$ exists, a contradiction.
If $f=2$, then $\Ga_N$ is isomorphic to one of the graphs in Lemma \ref{insoluble-2}, which also contradicts to Lemma \ref{Lemma2}.
This completes the proof of the lemma.
\qed



\end{document}